\documentclass{article}
\usepackage{a4wide} 
\usepackage{amssymb,amsmath,nccmath}
\usepackage{graphicx}

\newtheorem{theorem}{Theorem}
\newtheorem{proposition}{Proposition}
\newtheorem{lemma}{Lemma}
\newtheorem{coro}{Corollary}

\newtheorem{remark}{Remark}

\begin{document}

\begin{center}
\begin{Large}\begin{bf}
Hexagonal Inflation Tilings and Planar Monotiles\vspace{3ex}
\end{bf}\end{Large}

\begin{large}
Michael Baake$^{1}$, Franz G\"{a}hler$^{1}$ and Uwe Grimm$^{2}$\vspace{3ex}
\end{large}

\begin{small}
$^{1}$ Fakult\"{a}t f\"{u}r Mathematik, Universit\"{a}t Bielefeld,\\
       Postfach 100131, 33501 Bielefeld, Germany\\[0.5ex]
$^{2}$ Department of Mathematics and Statistics, The Open University,\\
       Walton Hall, Milton Keynes MK7 6AA, United Kingdom\vspace{1ex}
\end{small}
\end{center}

\begin{quote}
  Aperiodic tilings with a small number of prototiles are of
  particular interest, both theoretically and for applications in
  crystallography. In this direction, many people have tried to
  construct aperiodic tilings that are built from a single prototile
  with nearest neighbour matching rules, which is then called a
  monotile.  One strand of the search for a planar monotile has
  focussed on hexagonal analogues of Wang tiles. This led to two
  inflation tilings with interesting structural details. Both possess
  aperiodic local rules that define hulls with a model set
  structure. We review them in comparison, and clarify their relation
  with the classic half-hex tiling. In particular, we formulate
  various known results in a more comparative way, and augment them
  with some new results on the geometry and the topology of the
  underlying tiling spaces.
\end{quote}

\section{Introduction}

\begin{figure}[h]
\centerline{\includegraphics[width=0.4\textwidth]{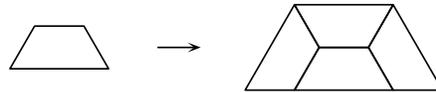}}
\caption{Half-hex inflation rule.}
\label{halfhex}
\end{figure}

A well-known inflation rule with integer inflation factor is the
half-hex inflation from \cite{GS}, Exercise~10.1.3 and Figure~10.1.7,
which we show in Figure~\ref{halfhex}. As such, it is a lattice
substitution (or inflation) in the sense of \cite{LMold,FS}. Moreover,
it is a face to face \emph{stone inflation} (in the sense of Danzer),
which means that each inflated tile is precisely dissected into copies
of the prototile so that the final tiling is face to face.  This rule
defines an aperiodic tiling of the plane, but it does not originate
from an aperiodic prototile set (for the terminology, we refer to
\cite{BGnew} and references therein). In principle, the procedure of
\cite{chaim} can be applied to add local information to the prototile
and to the inflation rule (via suitable markers and colours), until
one arrives at a version with an aperiodic prototile set. However, to
our knowledge, this has never been carried out, as it (most likely)
would result in a rather large prototile set.

Interestingly, two different inflation rules for hexagonally shaped
prototiles have independently been constructed, namely one by Roger
Penrose \cite{Pen97} and one by Joan Taylor \cite{Tay10} (see also
\cite{ST11}), each defining a tiling hull that can also be
characterised by aperiodic local rules. Viewed as dynamical systems
under the translation action of $\mathbb{R}^2$, they both possess the
continuous half-hex hull (and also the arrowed half-hex, to be
introduced later) as a topological factor, though with subtle
differences.  They may be considered as covers of the half-hex that
comprise just enough local information to admit an aperiodic prototile
set. In fact, both examples are again lattice inflations. They were
found in the attempt to construct an aperiodic planar monotile,
which loosely speaking is a single prototile together with some
local rules that tiles the plane, but only non-periodically. Let us 
mention that Taylor's original inflation tiling can be embedded into a
slightly larger tiling space that still possesses aperiodic local
rules \cite{ST11}, but is no longer minimal (see below for more on
minimality).

All three, the half-hex, the Penrose and the Taylor tiling, are
structures that can be described as model sets; see
\cite{M97,Moody00,TAO} for background on model sets. For the half-hex
tiling, this was first shown in \cite{Fre02}. Observing that two
half-hexes always join to form a regular hexagon, with edge length $1$
say, one obtains a hexagonal packing where three types of hexagons are
distinguished by a single diagonal line. We represent each hexagon by
a single point located at its centre, of type $\ell\in\{0,1,2\}$,
where $\ell$ corresponds to a diagonal that is rotated by $\ell\pi/3$
against the horizontal. This gives a partition
\[
    H^{}_{0}\cup H^{}_{1}\cup H^{}_{2}
    \, = \, \varGamma \, = \, 
    \sqrt{3}\, \langle\xi,\xi^{3}\rangle^{}_{\mathbb{Z}}
\]
with $\xi=e^{\pi i/6}$, where $\varGamma$ is a triangular lattice of
density $\frac{2}{9}\sqrt{3}$. When starting from a fixed point tiling
of the half-hex inflation rule, the fixed point equations for the
three point sets lead to the solution
\begin{equation}\label{halfhexsol}
   H^{}_{\ell} \, = \, A^{}_{\ell} \cup
   \bigcup_{n\ge 0} 2^{n} (2\varGamma + \sqrt{3}\,\xi^{3+2\ell})\, ,
\end{equation}
where the $A^{}_{\ell}$ are empty sets except for one, which is the
singleton set $\{0\}$. Where the latter occurs depends on the seed of
the selected fixed point. This means that there are precisely three
possibilities, corresponding to the three possible choices for the
central hexagon; see \cite{TAO} for the detailed derivation. 

This description establishes the model set structure, with the
$2$-adic completion of $\varGamma$ as internal space. Note that the
point $0$ is the unique limit point of any of the three unions in
Eq.~\eqref{halfhexsol} in the $2$-adic topology. The structure of the
union over expanded and shifted copies of $\varGamma$ is also called a
\emph{Toeplitz structure}. It is an example of a limit-periodic system
\cite{GK} with pure point diffraction (and, equivalently, with pure
point dynamical spectrum \cite{LMS,Martin,BL}). The diffraction
measure of a weighted Dirac comb on the half-hex tiling can be
calculated from Eq.~\eqref{halfhexsol} via the Poisson summation
formula, by an application of the methods explained in \cite{BM,TAO}.
This description of the half-hex tiling will be the key observation to
also identify the two covers (by Penrose and by Taylor) as model sets.

Below, we discuss the two inflation tilings due to Penrose and Taylor
in some detail. To be more precise with the latter case, we only
consider the minimal part of the tiling space considered by Socolar
and Taylor in \cite{ST11}. This minimal part is the tiling LI (local
indistinguishability) class defined by Taylor's original stone
inflation rule \cite{Tay10}. To distinguish the two hulls, we use the
term \emph{Taylor tiling} for the minimal inflation hull and refer to
the elements of the larger tiling space as the \emph{Socolar-Taylor
  tilings}.

We assume the reader to be familiar with the concept of \emph{mutual
  local derivability} (MLD), which was introduced in \cite{BSJ}; see
also \cite{B,TAO}. When the derivation rules commute with all symmetries
of the tilings under consideration (respectively their hulls), they
are called symmetry preserving. The corresponding equivalence classes
are called SMLD classes \cite{B}. 

The \emph{hull} of a (planar) tiling $\mathcal{T}$ with \emph{finite
  local complexity} (FLC) in $\mathbb{R}^{2}$ is defined as the orbit
closure in the local topology, $\mathbb{X} (\mathcal{T}) = \overline{
  \{ t+\mathcal{T} \mid t \in \mathbb{R}^{2} \} }$. Here, two tilings
are $\varepsilon$-close when they agree on the ball $B_{1/\varepsilon}
(0)$, possibly after (globally) translating one of them by an element
from $B_{\varepsilon} (0)$. Due to the FLC property, the hull is
compact \cite{Martin}, with continuous action of the group
$\mathbb{R}^{2}$ via translation. Consequently, the pair $(\mathbb{X}
(\mathcal{T}), \mathbb{R}^{2})$ is a topological dynamical system. It
is called \emph{minimal} when the translation orbit of every element
of the hull is dense in it. Our examples below will be minimal hulls
of FLC tilings, or of equivalent representatives of the corresponding
MLD class.

A hull $\mathbb{X}$ is called \emph{aperiodic} when no element of it
possesses non-trivial periods. In other words, $\mathbb{X}$ is
aperiodic when, for every $X\in\mathbb{X}$, the equation $t+X = X$
only holds for $t=0$. A hull is said to have \emph{local rules} when
it is specified by a finite list of legal local configurations, for
instance in the form of a finite atlas of patches.  If a set of local
rules specifies an aperiodic hull, the rules themselves are called
\emph{aperiodic}. When a set of rules specifies a hull that is
minimal, they are called \emph{perfect}. Of special interest now are
local rules that are aperiodic and perfect, such as the well-known
arrow matching conditions of the classic rhombic Penrose tiling;
compare \cite{GS}.

A single prototile, assumed compact and simply connected, is called a
\emph{monotile} (in the strict sense) when a set of aperiodic perfect
local rules exists that can be realised by nearest neighbour matchings
only. A subtle question in this context is whether one allows
reflected copies of the prototile or not.  Quite often, geometric
matching conditions are replaced by suitable decorations of the
prototile together with rules how these decorations have to form local
patterns in the tiling process.  The latter need not be restricted to
conditions for nearest neighbour tiles, in which case one speaks of a
\emph{functional monotile} to indicate the slightly more general
setting. Our planar examples below are of the latter type, or even a
further extension of it.  \medskip

This article, which is a brief review together with some new results
on the two tiling spaces, grew out of a meeting on discrete geometry
that was held at the Fields Institute in autumn 2011. As such, it is
primarily written for a readership with background in discrete
geometry, polytopes and tilings. We also try to provide the concepts
and methods for readers with a different background, though this is
often only possible by suitable pointers to the existing
literature. In particular, where correct and complete proofs are
available, we either refer to the original source or sketch how the
arguments have to be applied to suit our formulation. As is often the
case in discrete geometry, following a proof might need some pencil
and paper activity on the side of the reader; compare the introduction
and the type of presentation in \cite{GS}.  This is particularly true
of arguments around local derivation rules, inflation properties and
aperiodic prototile sets.

While we prepared this manuscript, we rewrote known results on both 
tiling spaces in a way that emphasises their similarities, and mildly
extended them, for instance by the percolation property of two derived
parity patterns. Moreover, we calculated several topological invariants
of the tiling spaces under consideration, which (as far as we are aware)
were not known before. In order not to create an imbalance, we only
describe how to do that in principle (again with proper references)
and then state the results. We also include the dynamical zeta functions
for the inflation action on the hulls, which are the generating 
functions for the corresponding fixed point counts. They turn out
to be particularly useful for deriving the structure of the hull.
Since several examples of a similar nature have recently been 
investigated in full detail, compare \cite{BGG-gtm,BGG-All,G-Rob}, 
we felt that this short account is adequate (in particular, as the 
explicit results are the outcome of a computer algebra program).
\medskip

The paper is organised as follows. In Section~\ref{sec:pen}, we begin
with a brief discussion of the\break $(1+\varepsilon + \varepsilon^{2})$-tiling
due to Penrose, together with various other elements of the MLD class
defined by it. Section~\ref{sec:Taylor} contains the corresponding
material on the Taylor (and the Socolar-Taylor) tilings and their
`derivatives', which we discuss in slightly more detail, including the
percolation result on the parity patterns of both tilings. The
topological invariants and various other quantities for a comparison
of the tilings are presented in Section~\ref{sec:inv}, which is
followed by some concluding remarks and open problems.

\section{Penrose's aperiodic hexagon tiling
and related patterns}\label{sec:pen}

\begin{figure}[t]
\centerline{\includegraphics[width=0.6\textwidth]{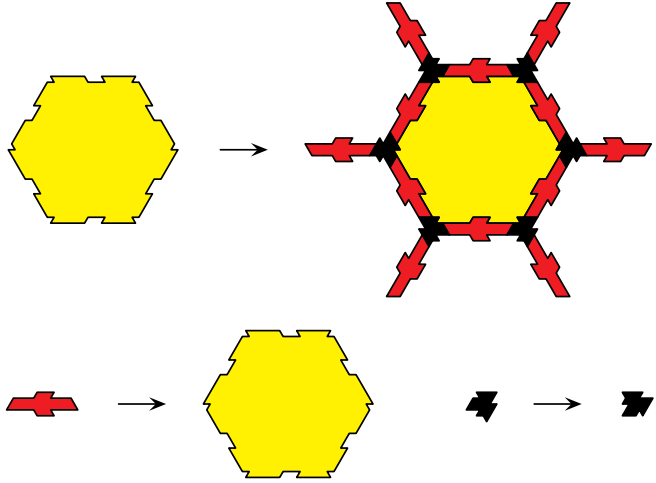}}
\caption{Inflation rule for Penrose's
  $(1+\varepsilon+\varepsilon^{2})$-tiling.}
\label{penepssub}
\end{figure}

\begin{figure}[t]
\centerline{\includegraphics[width=0.9\textwidth]{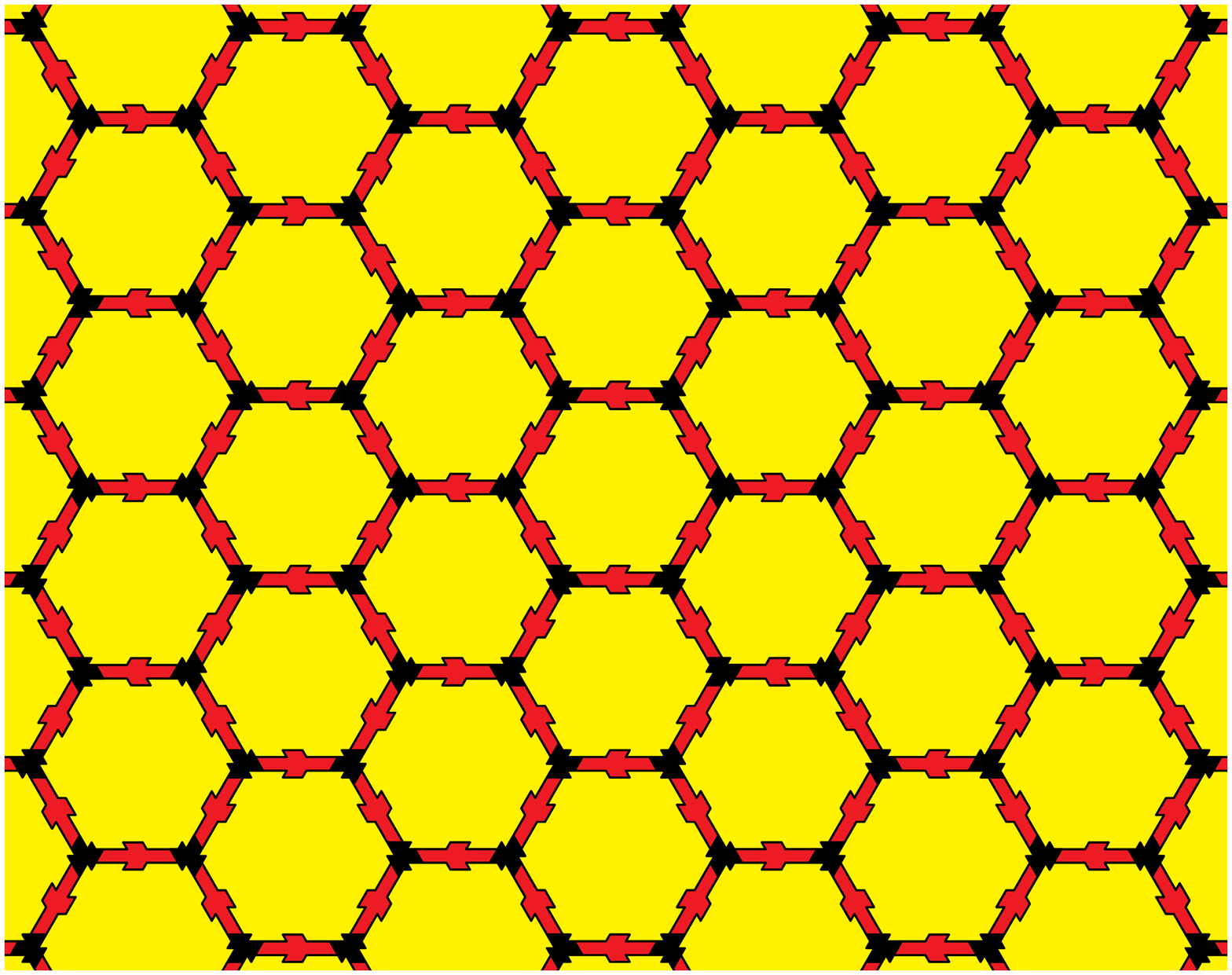}}
\caption{Patch of  Penrose's
  $(1+\varepsilon+\varepsilon^{2})$-tiling.}
\label{penepstil}
\end{figure}

The $(1+\varepsilon+\varepsilon^{2})$-tiling due to Roger Penrose is
built from three prototiles, up to Euclidean motions (including
reflections). The tiles and the inflation rule (with linear inflation
multiplier $2$) are shown in Figure~\ref{penepssub}. The name
$(1+\varepsilon+ \varepsilon^{2})$-tiling refers to the three
prototiles as the $1$-tile (hexagon), the $\varepsilon$-tile (edge
tile) and the $\varepsilon^{2}$-tile (corner tile). The
$\varepsilon$-tile has a definite length, but can be made arbitrarily
thin, while the $\varepsilon^{2}$-tile can be made arbitrarily small.
This inflation rule is not primitive, but still defines a unique
tiling LI class in the plane via a fixed point tiling with a hexagon
at its centre. A patch of such a tiling is shown in
Figure~\ref{penepstil}. There are $12$ fixed point tilings of this
type, each defining the same LI class.  These fixed points form a
single orbit under the $D_6$ symmetry of the LI class.

The following result was shown in \cite{Pen97} by the
composition-decomposition method; see \cite{G} for a detailed
description of this method.

\begin{proposition}\label{prop-one}
  The inflation rule of Figure~$\ref{penepssub}$ defines a unique
  tiling LI class with perfect aperiodic local rules. The latter are
  formulated via an aperiodic prototile set, which consists of the
  three tiles from Figure~$\ref{penepssub}$ $($together with rotated
  and reflected copies$)$.  The rules are then realised as purely
  geometric matching conditions of the tiles.  \hfill$\square$
\end{proposition}

\begin{figure}[t]
\centerline{\includegraphics[width=0.9\textwidth]{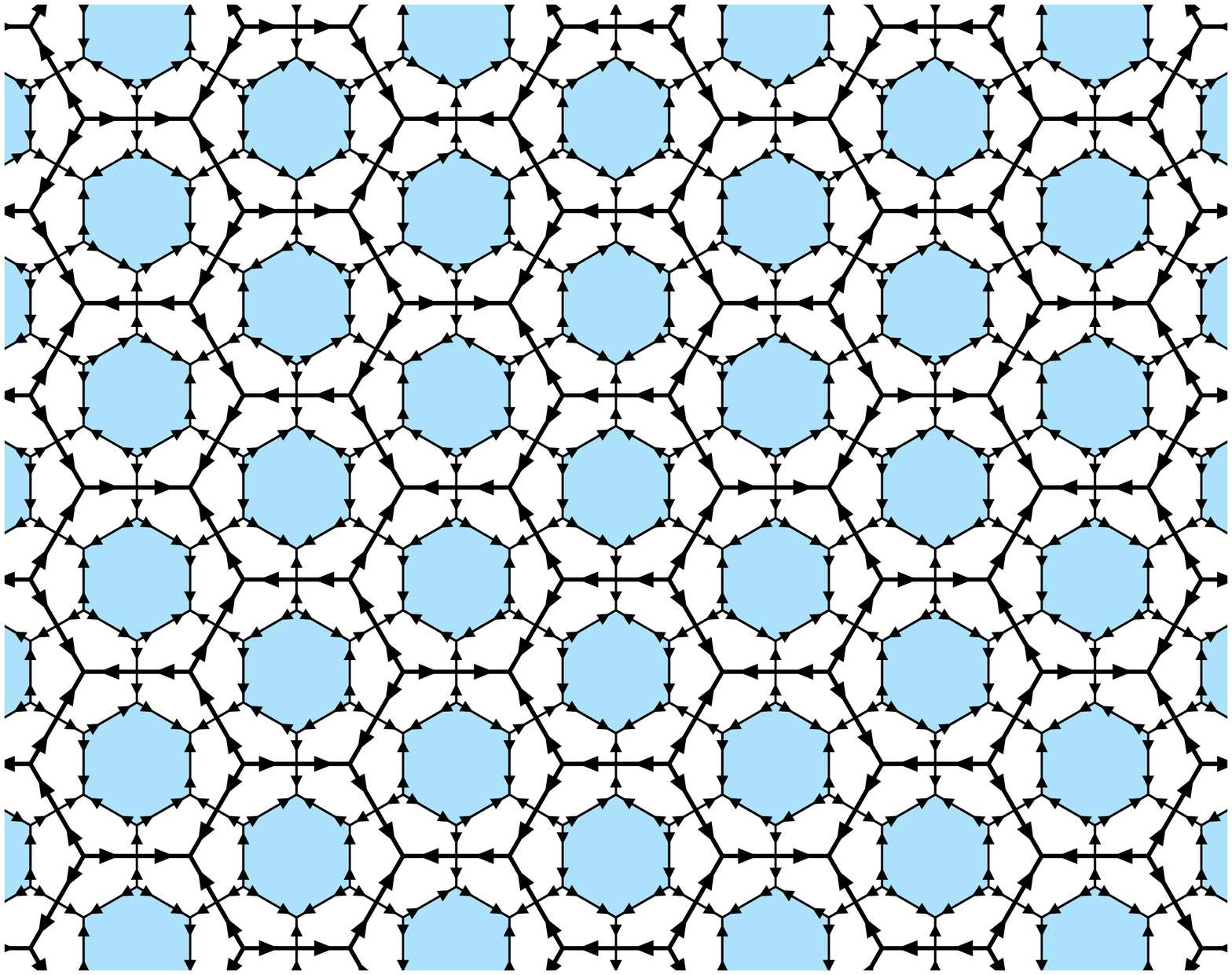}}
\caption{A patch of the double hexagon tiling, exactly corresponding 
to that of Figure~\ref{penepstil}.}
\label{penmono}
\end{figure}

\begin{figure}[b]
\centerline{\includegraphics[width=0.9\textwidth]{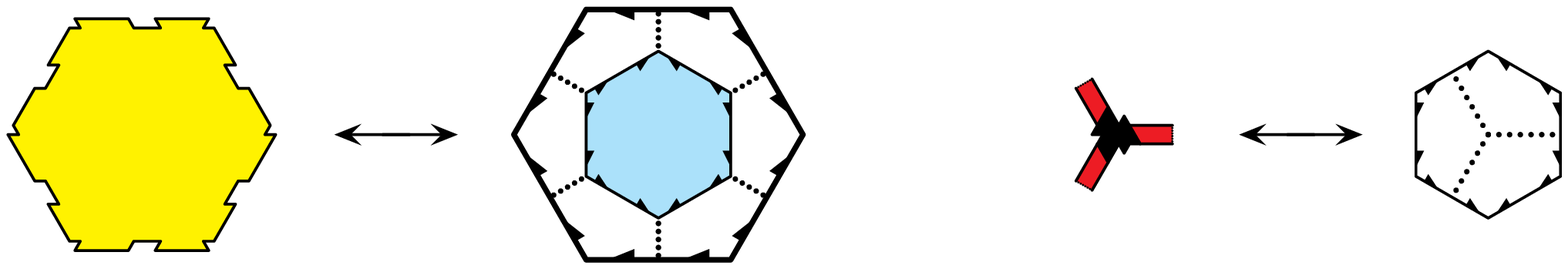}}
\caption{Rules for the mutual local derivation between 
  $(1+\varepsilon+\varepsilon^{2})$-tilings and double hexagon
  tilings.}
\label{penmonomld}
\end{figure}

In the original publication \cite{Pen97}, it was argued that this
system comes `close' to an aperiodic monotile in the sense that it
essentially is a marked hexagon tiling with matching conditions that
are realised by `key tiles', which can be made thin and small.  A
transformation to an equivalent version was only sketched briefly at
the end of the article, and subsequently substantiated in
\cite{PenTwist} in form of a puzzle and its solution. The key idea is
to change from the $(1+\varepsilon+\varepsilon^{2})$-tiling to the
double hexagon tiling of Figure~\ref{penmono}, which is possible by
the local derivation rule sketched in Figure~\ref{penmonomld}, when
read from left to right.

The actual rule consists of two steps.  In the first, each hexagonal
$1$-tile is replaced by a double hexagon as shown in the left panel of
Figure~\ref{penmonomld}. This produces a complete tiling of the larger
hexagons (with matching arrows) with fully oriented inscribed
hexagons, while all remaining small hexagons, namely those around the
vertices of the larger hexagons, are still incomplete (thin dotted
lines). Their missing orientations derive consistently from the
$\varepsilon^{2}$-tiles via the rule shown in the right panel.

Note that all hexagons (which means on both scales) have the same type
of arrow pattern. In particular, each hexagon has precisely one pair
of parallel edges that are oriented in the same direction, while all
other arrows point towards the two remaining (antipodal)
vertices. Note that this is the same edge orientation pattern as seen
at the boundary of the hexagonal $1$-tile in Figures~\ref{penepssub}
and \ref{penepstil}.

The derivation rule specified in Figure~\ref{penmonomld} is clearly
local in both directions, and commutes with the translation action as
well as with all symmetry operations of the group $D_{6}$ (the
symmetry group of the regular hexagon). The following result is thus
obvious.

\begin{proposition}
  The LI classes of the $(1+\varepsilon+\varepsilon^{2})$-tiling of
  Figure~$\ref{penepstil}$ and the Penrose double hexagon pattern of
  Figure~$\ref{penmono}$ are SMLD.\hfill$\square$
\end{proposition}

\begin{figure}[t]
\centerline{\includegraphics[width=0.9\textwidth]{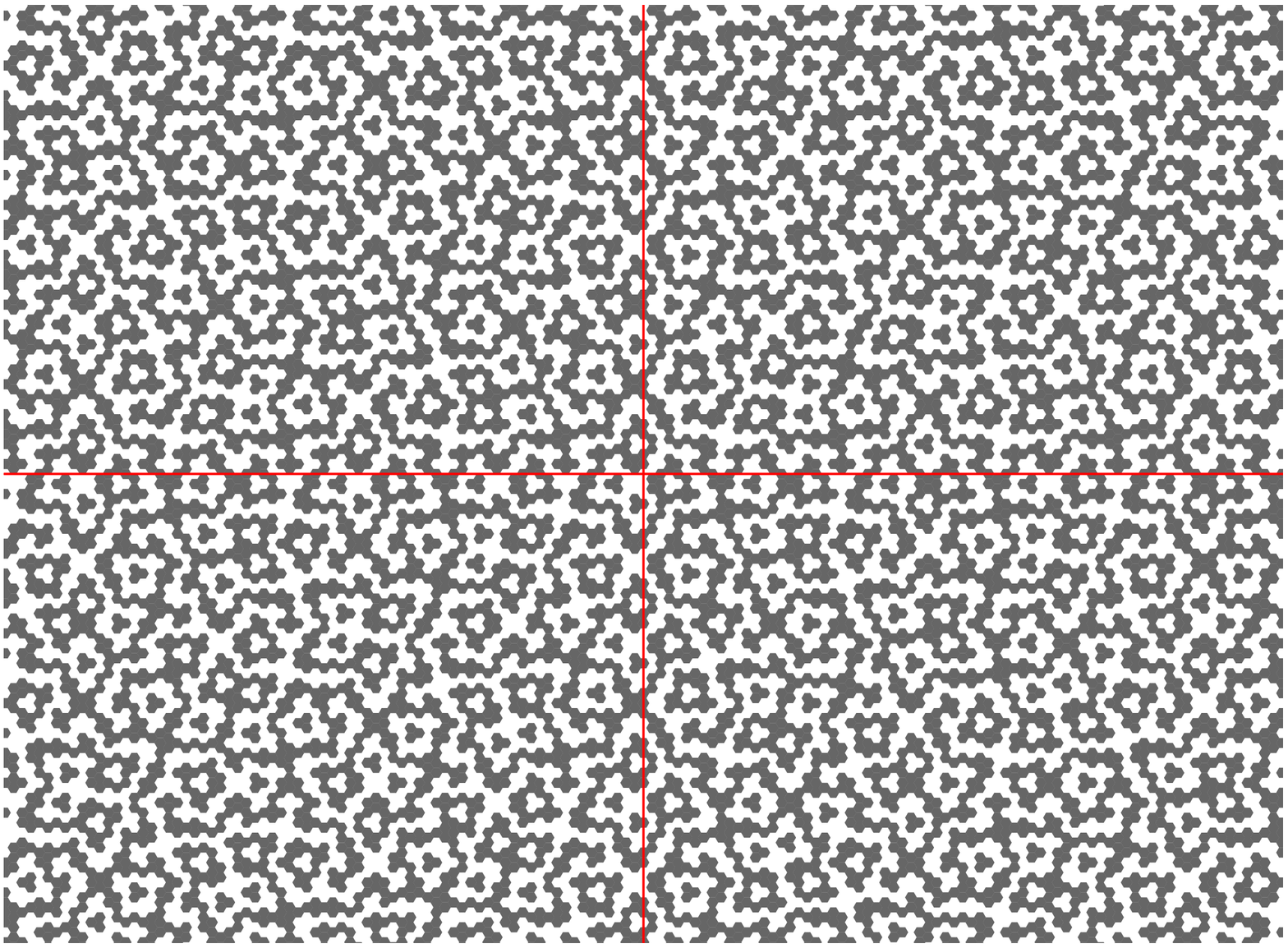}}
\caption{Parity pattern of a $(1+\varepsilon+\varepsilon^{2})$-tiling,
  as derived from a fixed point tiling of the inflation rule of
  Figure~\ref{penepssub}.  This particular pattern (also the infinite
  one) possesses an almost colour reflection symmetry with respect to
  the indicated lines; see text for details.}
\label{penepspar}
\end{figure}

A closer inspection of the prototile set of Figure~\ref{penepssub}
reveals that the hexagonal $1$-tile occurs in two chiralities with six
orientations each. Keeping track of the chiralities only (by two
colours, white and grey say) and disregarding all other structural
elements provides a local derivation from the $(1+\varepsilon +
\varepsilon^{2})$-tiling to an ensemble of $2$-colourings of the
hexagonal packing. An example of the latter is illustrated in
Figure~\ref{penepspar}. We call the elements of this new ensemble the
\emph{parity patterns} of the $(1+\varepsilon +
\varepsilon^{2})$-tilings. As far as we are aware, this parity pattern
has not been considered before. Our motive to do so will become
clear from the comparison with the llama tilings in the next section.
By construction, the parity patterns form a single LI
class. A little surprising is the following property.

\begin{theorem}\label{thm1}
  The LI class of the $(1+\varepsilon+\varepsilon^{2})$-tiling of
  Figure~$\ref{penepstil}$ and that of the corresponding parity patterns
  of Figure~$\ref{penepspar}$ are SMLD.
\end{theorem}
\noindent\textsc{Sketch of Proof.}
The determination of the parity pattern that belongs to
a $(1+\varepsilon+\varepsilon^{2})$-tiling is clearly local
and preserves all symmetries, so that this direction is clear.

Conversely, starting from a parity pattern, the corresponding
$(1+\varepsilon+\varepsilon^{2})$-tiling is locally reconstructed via
hexagonal coronae of order $3$. A simple computer search shows that
such coronae uniquely specify the decorated hexagon that corresponds
to its centre, once again in a symmetry-preserving way.
\hfill$\square$\bigskip

Let us note that, as a result of the relation between two fixed point
tilings with mirror image central hexagons, the particular parity
patch of Figure~\ref{penepspar} shows an almost reflection colour
symmetry for the reflections in the two lines indicated.  More
precisely, under reflection and colour inversion, the patch is mapped
onto itself, except for some hexagons along the reflection line.

\begin{remark}
  The left panel of Figure~$\ref{penmonomld}$ shows the building block
  of the double hexagon tiling without arrows on the dashed lines.  As
  mentioned above, they are added by the local rule how to complete
  the vertex configurations. In this sense, one has a single prototile
  \emph{template} together with a set of local rules that specify how
  to put them together. Any double hexagon tiling of the plane that
  everywhere satisfies the rules is an element of the double hexagon
  LI class, and in this sense one has a functional monotile
  template. We suggest calling this a \emph{weak functional monotile},
  as it stretches the monotile concept to some extent.
\end{remark}

Let us return to the structure of the $(1+\varepsilon +
\varepsilon^{2})$-tilings. As mentioned before, each hexagonal
$1$-tile possesses a unique pair of parallel edges that are
oriented (by the arrows) in the same direction. If we divide it
now into two half-hexes along the diagonal that is parallel to
this edge pair, we can locally derive a half-hex tiling from any
$(1+\varepsilon + \varepsilon^{2})$-tiling. An inspection of
Figure~\ref{penepssub} confirms that this derivation rule indeed
induces the half-hex inflation of Figure~\ref{halfhex}. This is
a local derivation of sliding block map type on the underlying
hexagonal packing, hence continuous in the local topology.

\begin{proposition}
  The LI class of the half-hex tiling as defined by the stone
  inflation rule of Figure~$\ref{halfhex}$ defines a minimal
  topological dynamical system under the translation action of\/
  $\mathbb{R}^{2}$ that is a topological factor of the LI class of the
  $(1+\varepsilon + \varepsilon^{2})$-tilings defined by the inflation
  rule of Figure~$\ref{penepssub}$. The corresponding factor map is
  one-to-one almost everywhere, but the two tiling spaces define
  distinct MLD classes.
\end{proposition}
\noindent\textsc{Sketch of Proof.}
The first claim follows from our above description of the local
derivation rule via the additional diagonal line in the hexagonal
$1$-tiles. The two LI classes cannot be MLD (and hence also not SMLD)
because the existence of aperiodic local rules is an invariant
property of an MLD class, hence shared by all LI classes that are
MLD. It is well-known \cite{GS,Fre02} that the half-hex hull has no such
set of rules, because any finite atlas of patches can still be part of
a periodic arrangement.

Since the $(1+\varepsilon + \varepsilon^{2})$-tilings have perfect
aperiodic local rules by Proposition~\ref{prop-one}, the last claim is
clear. The statement on the multiplicity of the mapping is a
consequence of the model set structure, which we prove below in
Theorem~\ref{pure-point}.
\hfill$\square$\bigskip

The (regular) model set structure of the half-hex tiling, as spelled
out in Eq.~\eqref{halfhexsol} for the fixed points under the inflation
rule, implies that there is a `torus parametrisation' map onto a
compact Abelian group \cite{Martin,BLM}. Here, it is a factor map onto
the two-dimensional dyadic solenoid $\mathbb{S}^{2}_{2}$, which is
almost everywhere one-to-one by \cite[Thm.~5]{BLM}. Since the half-hex
LI class is the image of the $(1+\varepsilon + \varepsilon^{2})$ LI
class under a factor map that is itself one-to-one almost everywhere
(which follows by the same argument that we use below to prove
Theorem~\ref{pure-point}), we know (via concatenation) that there
exists an almost everywhere one-to-one factor map from the
$(1+\varepsilon + \varepsilon^{2})$ LI class onto
$\mathbb{S}^{2}_{2}$. Then, \cite[Thm.~6]{BLM} implies the following
result.

\begin{coro}
  The\/ $(1+\varepsilon + \varepsilon^{2})$ LI class has a model set
  structure, with the same cut and project scheme as derived for the
  half-hex tilings. \hfill$\square$
\end{coro}

In summary, the $(1+\varepsilon + \varepsilon^{2})$ LI class has all
magical properties: It can be defined by an inflation rule, by a
set of perfect aperiodic local rules, and as a regular model set.
Moreover, it comes close to solving the quest for a (functional)
monotile.

Let us turn our attention to a later (though completely independent)
attempt of a similar kind that improves the monotile state-of-affairs.

\section{Taylor's inflation tiling}\label{sec:Taylor}

\begin{figure}[t]
\centerline{\includegraphics[width=0.65\textwidth]{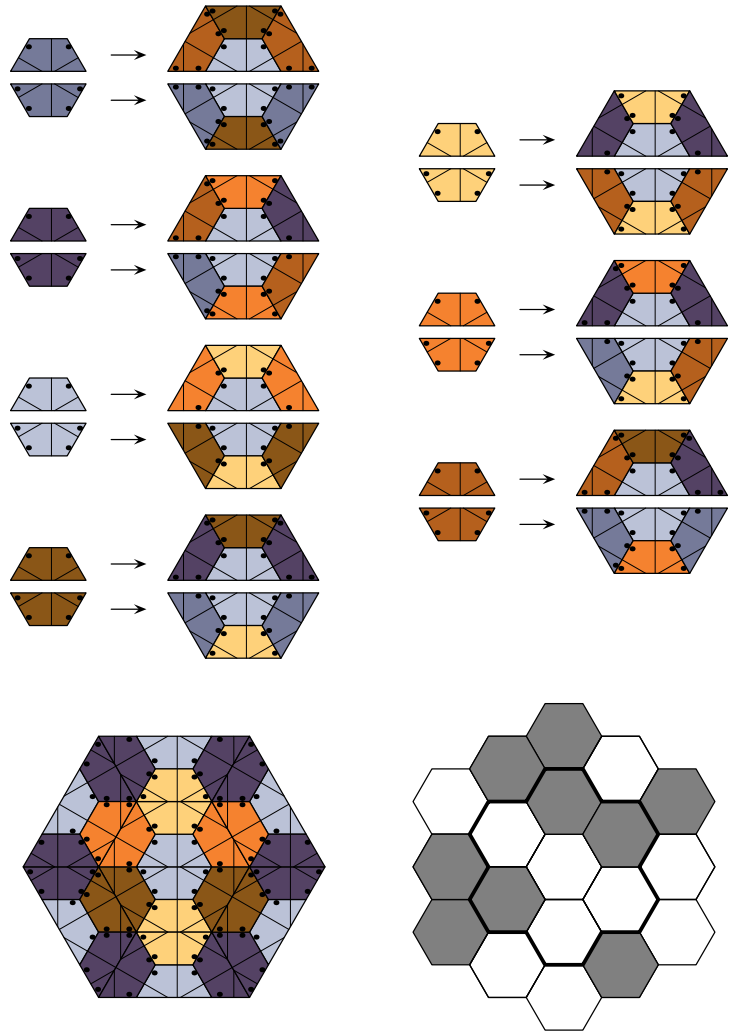}}
\caption{The primitive inflation rule of Taylor's half-hex 
inflation (top). The central patch of a fixed point tiling is shown
on the lower left panel, with its parity pattern to the right.}
\label{taysub}
\end{figure}

Consider the primitive stone inflation rule of Figure~\ref{taysub}.
It is formulated with $14$ prototiles (up to similarity) of half-hex
shape that are distinguished by colour and a decoration (with points
and lines). Each half-hex occurs in two chiralities and six
orientations, so that the total number of prototiles (up to
translation) is $168$. In the original paper \cite{Tay10}, the $7$
colours are labelled $A,B,C, \ldots, G$. The $C$-type tiles (light blue
in our version) are special in the sense that they are in the centre
of any fixed point tiling under the inflation rule. They are also more
frequent than the other types. As mentioned above, we call the
elements of the tiling space (or hull) defined by this inflation rule
the \emph{Taylor tilings}. The following result is immediate from
Figure~\ref{taysub}, in comparison with Figure~\ref{halfhex}, by
simply removing decorations and colour.

\begin{lemma}
   The half-hex LI class is a topological factor of the LI class
   of the Taylor tilings, where the latter is again minimal.
   In particular, the Taylor LI class is aperiodic. 
   \hfill $\square$
\end{lemma}

\begin{figure}[t]
\centerline{\includegraphics[width=0.8\textwidth]{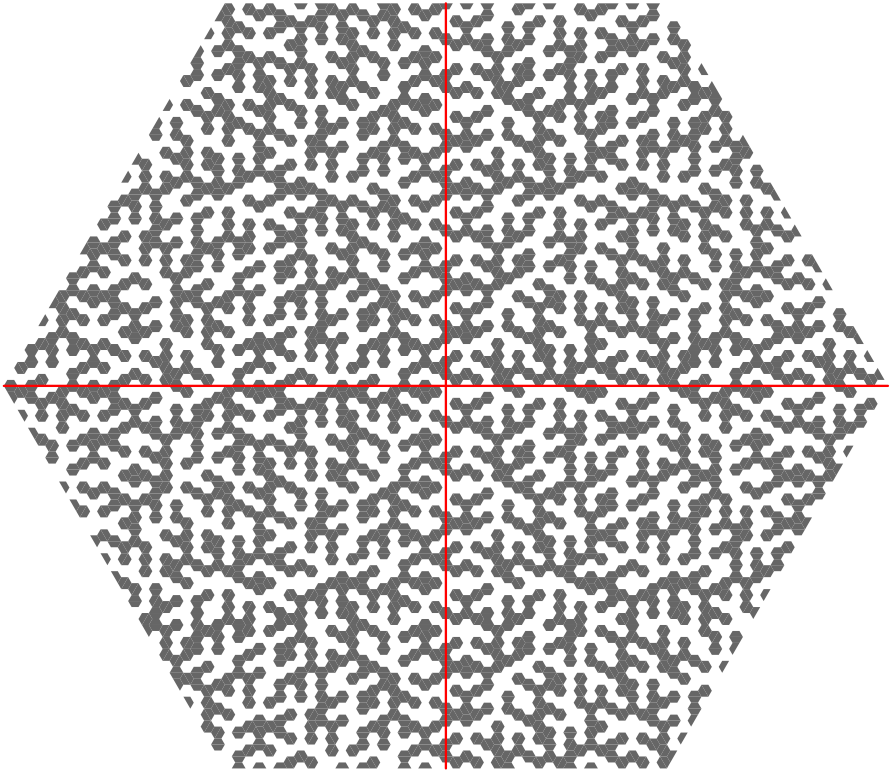}}
\caption{Patch of a llama tiling, as derived from a fixed point tiling
  of the inflation rule of Figure~\ref{taysub}. This particular
  pattern (and its infinite extension) possesses an almost colour
  reflection symmetry with respect to the indicated lines; see text
  for details.}
\label{llama}
\end{figure}

As in the case of the Penrose tiling above, one can locally derive a
hexagonal parity pattern from every Taylor tiling. To this end, one
considers the natural half-hex pairs, disregards their colours, and
applies a grey/white coding of the two chiralities. The resulting
two-coloured hexagonal packings are called \emph{llama tilings}, see
Figure~\ref{llama} for an illustration. The name refers to the shape
of the smallest island (of either colour), and was coined by
Taylor. As in the previous section, the parity pattern still contains
the full local information.  Also, the llama tiling has the same type
of almost colour reflection symmetry that we encountered in the
previous section for the Penrose parity pattern.

\begin{theorem}
  The LI class of Taylor's inflation tiling and that of the llama
  tiling of Figure~$\ref{llama}$ are SMLD. In particular, the llama
  tiling is aperiodic.
\end{theorem}
\noindent\textsc{Sketch of Proof.}
As in the case of Theorem~\ref{thm1}, the derivation of the llama
tiling, which is the parity pattern of the Taylor tiling, is obviously
local and symmetry preserving. 

\begin{figure}[t]
\centerline{\includegraphics[width=0.8\textwidth]{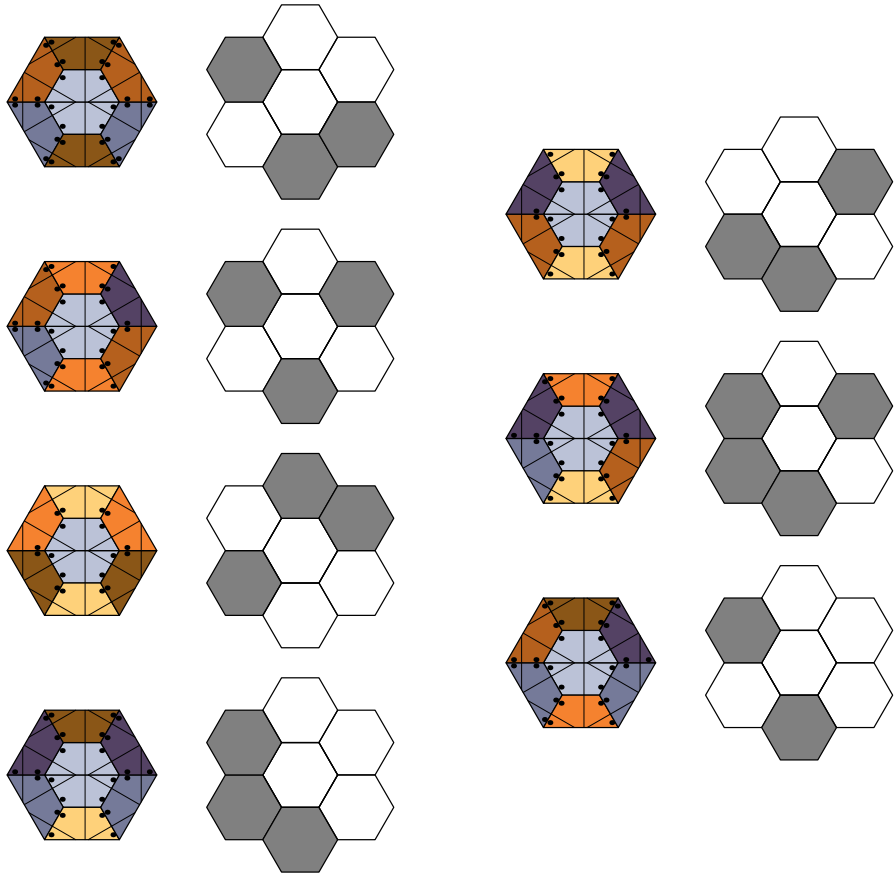}}
\caption{Correspondence between fully decorated (half) hexagons,
with central hexagon of type $C$, and their parity patterns.}
\label{tayrec}
\end{figure}

For the converse direction, there are three proofs known. The first is
based on an idea by Joan Taylor, and is spelled out in detail in
\cite{LM}. A related argument uses the local information contained in
the llamas together with the correspondence of coloured and marked
hexagons with local parity patterns, as indicated in
Figure~\ref{tayrec}; see \cite{TAO} for details. Finally, as for
Theorem~\ref{thm1}, one can reconstruct the complete decoration of any
hexagon from the order-$3$ coronae of the llama tiling, which
preserves the symmetry.

The final claim follows from the aperiodicity of the Taylor tiling 
by standard arguments.
\hfill$\square$\bigskip

Before we continue with our general discussion, let us mention
an interesting property of the llama tilings, which also holds
for the Penrose parity pattern.

\begin{theorem}\label{taypatches}
  The llama tiling of Figure~$\ref{llama}$ possesses connected
  components of either colour of unbounded size. The same 
  conclusion also holds for the Penrose parity pattern of
  Figure~$\ref{penepspar}$. Moreover, this property extends
  to every element of the respective LI classes.
\end{theorem}
\noindent\textsc{Proof.}
The patch of the Taylor tiling shown in Figure~\ref{llama} is derived
from one inflation fixed point with a hexagon (originally of type $C$)
of positive chirality as a seed, which we denote as pattern
P1. Another fixed point pattern, P2, can be obtained from the
corresponding hexagon of opposite chirality, which has grey and white
colours interchanged relative to P1. Nevertheless, P1 and P2 are LI,
so that arbitrarily large patches of either pattern occur in the
other.

Now, assume that P1 does not contain connected patches of white
hexagons of unbounded size, where two hexagons are called connected
when they share an edge. If so, there must be a maximal white `island'
in P1, of diameter $r$ say, which must then be surrounded by a
connected `belt' of grey hexagons (which is possibly part of an even
larger connected patch). Then, the same belt exists in P2, this time
as a white belt around a grey island, and this belt has diameter $>r$
by construction. Since P1 and P2 are LI, this patch from P2 must also
occur somewhere in P1, in contradiction to the assumption, and our
claim follows.

The argument for the Penrose parity pattern is completely analogous
to that for the llama tilings, while the final claim is obvious.
\hfill$\square$\bigskip

The llama tiling as well as Penrose's parity pattern are thus
interesting examples of deterministic aperiodic structures with
\emph{percolation}. In particular, an element with an infinite
connected cluster must exist in either LI class.  Since each class is
a compact space, this claim follows from a compactness argument,
because any sequence of tilings with connected clusters of increasing
diameters around the origin must contain a subsequence that converges
in the local topology to a tiling with an infinite cluster. Note,
however, that the above argument does not imply the existence of a
sequence of islands of growing size, which has been conjectured for
the llama LI class \cite{ST11} and seems equally likely in the other
LI class as well, both on the basis of inflation series of suitable
patches.

\begin{figure}[t]
\centerline{\includegraphics[width=0.8\textwidth]{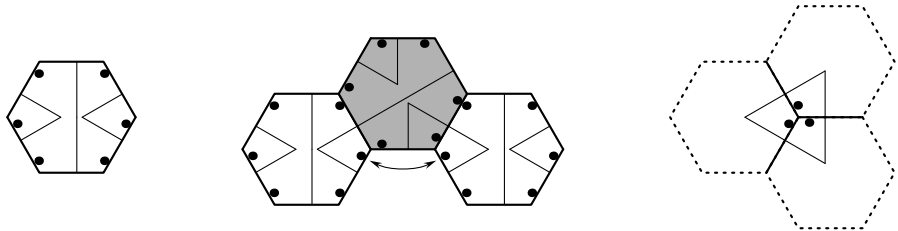}}
\caption{Taylor's functional monotile (left), a sketch of the edge
  transfer rule R2 (centre) and the forbidden threefold vertex seed
  (right); see text for details.}
\label{taymatch}
\end{figure}

Between any Taylor tiling and the corresponding llama tiling (which
are SMLD) is another version that still shows all line and point
markings (and hence the chirality of the hexagons), but not the seven
colours.  Clearly, also this version, which we call the \emph{decorated
llama tiling}, is in the same SMLD class. One can now formulate three
local rules for the corresponding prototile set, which consists of
$12$ tiles (up to translations).
\begin{itemize}
   \item[R1.] The hexagons must match at common edges in the sense that
     the decoration lines do not jump on crossing the common edge.
   \item[R2.] The point markers must satisfy the edge transfer rule
     sketched in the middle panel of Figure~\ref{taymatch}, as
     indicated by the arrow, for any pair of hexagons separated by a
     single edge. The two points at the corners adjacent to that edge
     have to be in the same position, and this rule applies
     irrespective of the chirality types of the tiles.
   \item[R3.] No vertex configuration is allowed to have adjacent
     points in a threefold symmetric arrangement, such as the one
     shown in the right panel of Figure~\ref{taymatch}, or its rotated
     and reflected versions.
\end{itemize}
In \cite{Tay10,ST11}, the following result is shown by the
composition-decomposition method, with immediate consequences
for other members of the corresponding MLD class.

\begin{theorem}
  The rules R1--R3 constitute perfect aperiodic local rules
  for the LI class of the decorated llama tiling. Consequently, a
  corresponding set of perfect aperiodic local rules also exists for
  the LI class of the Taylor tilings and for the llama LI class.
  \hfill $\square$
\end{theorem}

The decorated llama tiling was selected in the MLD class for the
following reason.

\begin{coro}
  The hexagon of the left panel of Figure~$\ref{taymatch}$, together
  with its reflected copy, provides a functional monotile under the
  local rules R1--R3. It defines the LI class of the
  decorated llama tiling.  \hfill $\square$
\end{coro}

Let us mention that one can relax the three rules by omitting R3.  As
is proved in \cite{ST11}, this still provides a set of aperiodic local
rules, this time for the larger space of the Socolar-Taylor
tilings. This space is not minimal, and additionally contains two
patterns with global threefold rotation symmetry (and their
translates).

\section{Topological invariants and the structure of the hulls}
\label{sec:inv}

In this section, we derive and compare further details of the hulls of
the various hexagonal tilings, where we go beyond the material that
has appeared in the literature so far.  For this purpose, it proves
useful to generate the tilings by lattice inflations in which the
tiles are represented by the points of a triangular lattice, with the
tile type attached as a label to each point. Geometrically, such
lattice inflation tilings consist of labelled hexagon tiles only.  Key
tiles like the ones in the $(1+\varepsilon+\varepsilon^2)$-tiling have
to be absorbed into the type of the tiles (via suitable decorations).

\begin{figure}[b]
\centerline{
\includegraphics[width=0.28\textwidth]{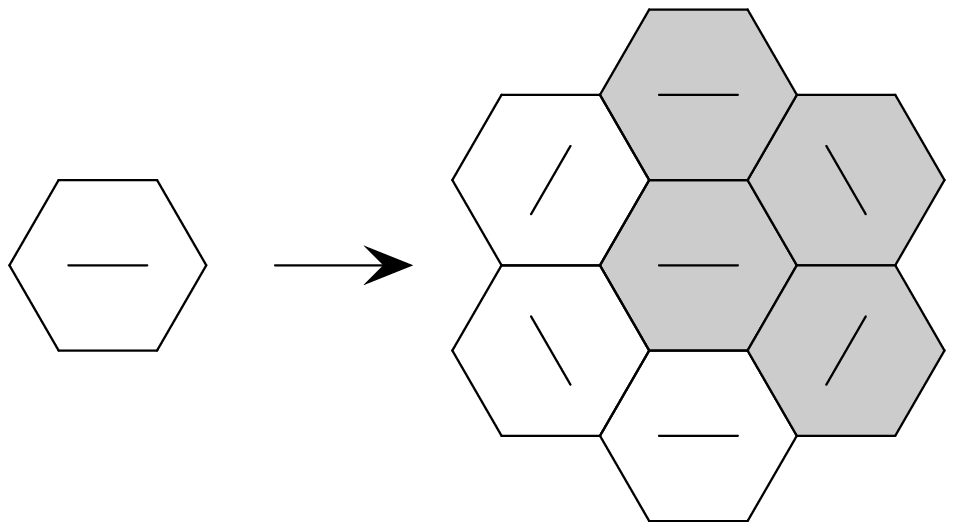}\hfill
\includegraphics[width=0.28\textwidth]{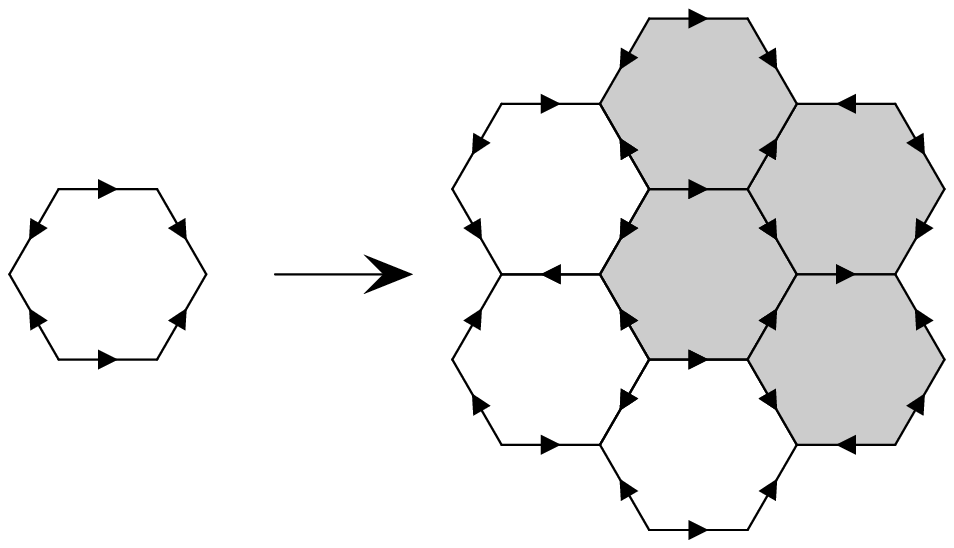}\hfill
\includegraphics[width=0.28\textwidth]{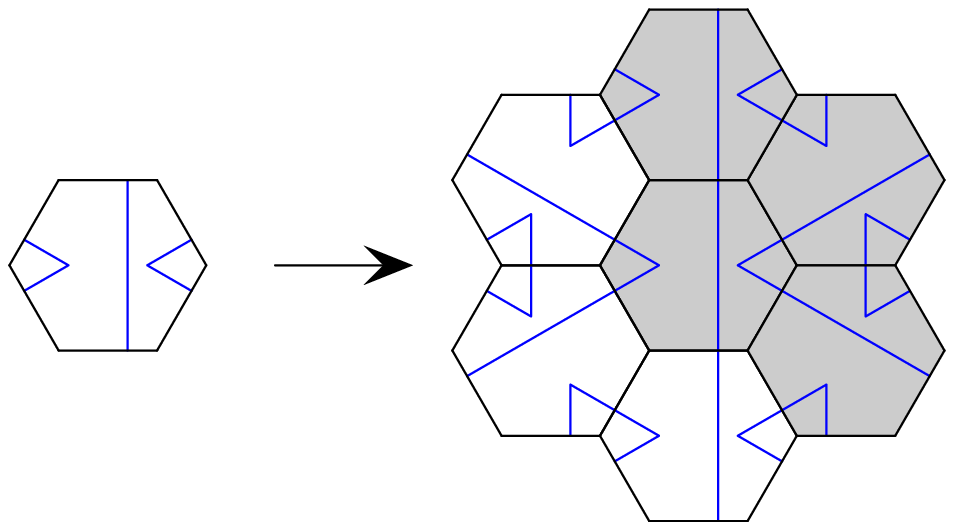}
}
\caption{Hexagon (pseudo) inflation for the half-hex (left) and for
  two (locally equivalent) variants of the arrowed half-hex (centre
  and right). If only the shaded hexagons are retained, a standard
  inflation for a wedge with opening angle $\frac{2\pi}{3}$ is obtained; see
  text for details.}
\label{hex-pseudo-inf}
\end{figure}

When passing from half hexagons to full hexagons, one easily obtains
the overlapping inflation (also called pseudo inflation) shown in
Figure~\ref{hex-pseudo-inf} (left) for the half-hex tiling. The outer
ring of hexagons is shared between the inflations of neighbouring
hexagons. A pseudo inflation has to be \emph{consistent} in the sense
that it agrees on the overlap regions. From this pseudo inflation, a
standard inflation can be obtained by replacing each hexagon by the
four shaded hexagons only, independently of the orientation of the
original hexagon. Conversely, each hexagon in the inflated tiling is
assigned to a unique supertile hexagon. The inflation rule obtained in
this way is not rotation covariant, but this is only a minor
disadvantage.  Also, the fixed point tiling obtained under the
iterated inflations no longer covers the entire plane, but only a
120-degree wedge, which is not a problem either. Since the
inflation rule is derived from a pseudo inflation, it determines the
tiling also in the other sectors, and there is a unique continuation
to the rest of the plane. Technically speaking, the inflation rule
\emph{forces the border} \cite{AP98}, a property which simplifies the
computation of topological invariants considerably.

A natural generalisation of the half-hex tiling is the \emph{arrowed
  half-hex tiling}, obtained by either of the two variants of (pseudo)
inflations shown in Figure~\ref{hex-pseudo-inf} (centre and
right). Clearly, the two decorations are locally equivalent, and lead
to tilings which are MLD.  Just as for the half-hex, a pseudo inflation
is obtained first, from which a standard inflation is then derived by
taking only the shaded tiles. The arrow pattern on the hexagons shown
in Figure~\ref{hex-pseudo-inf} (centre) is the same as that of the
hexagon tile of the $(1+\varepsilon+\varepsilon^2)$-tiling. In fact,
it is easy to see that the arrowed half-hex is locally derivable from
the latter, and thus is a factor. Similarly, the hexagon decoration
of Figure~\ref{hex-pseudo-inf} (right) is part of the decoration of the
Taylor hexagons, so that the arrowed half-hex is also a factor of the
Taylor tiling.

\begin{coro}
  The hull of the half-hex tiling, viewed as a dynamical system under
  the action of\/ $\mathbb{R}^{2}$, is a topological factor of both
  the $(1+\varepsilon+\varepsilon^2)$-tiling hull and the Taylor
  tiling hull.  The same property holds for the hull of the arrowed
  half-hex tiling.  \hfill $\square$
\end{coro}

\begin{figure}[t]
\centerline{
\includegraphics[width=0.6\textwidth]{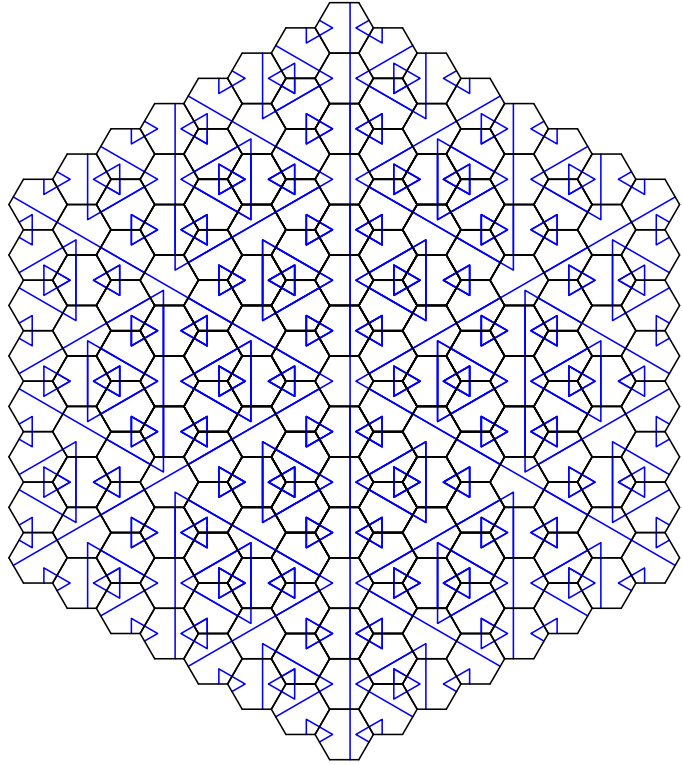}
}
\caption{Patch of the arrowed half-hex tiling. The blue decoration
lines form an infinite hierarchy of triangles of all sizes. The
blue triangles of a given size and orientation form a triangular
lattice.}
\label{arr-half-hex}
\end{figure}

A larger patch of the arrowed half-hex tiling is shown in
Figure~\ref{arr-half-hex} (in the Taylor decoration variant). The
fixed point tilings generated from the six different orientations of
the arrowed hexagon differ only along three rows of hexagons crossing
in the centre of the figure. We call these the \emph{singular
  rows}. The six fixed point tilings project to the same point on the
solenoid $\mathbb{S}^2_2$. Moreover, pushing the fixed point centre along one
of the singular rows towards infinity, pairs of tilings are obtained
which are mirror images of each other, except on the one singular row
that remains. Each such pair of tilings also projects to the same
points of the solenoid $\mathbb{S}^2_2$. In fact, each singular row
gives rise to a 1d sub-solenoid $\mathbb{S}^{1}_{2}$, onto which the
projection is 2-to-1 (except for the crossing point of these
sub-solenoids). Three such sub-solenoids $\mathbb{S}^{1}_{2}$, onto
which the projection is 2-to-1, clearly exist in the $(1+\varepsilon+
\varepsilon^2)$-tiling and in the Taylor tiling as well. However, for
the (naked) half-hex tiling, there are \emph{no} singular rows. In this
case, there are just three fixed point tilings, differing in the
orientation of the central hexagon only, which project to the same
point on the solenoid $\mathbb{S}^2_2$. We refer to \cite{TAO} for an
explicit derivation of the model set coordinates of the arrowed
half-hex tiling.

\begin{figure}[t]
\centerline{\includegraphics[width=0.6\textwidth]{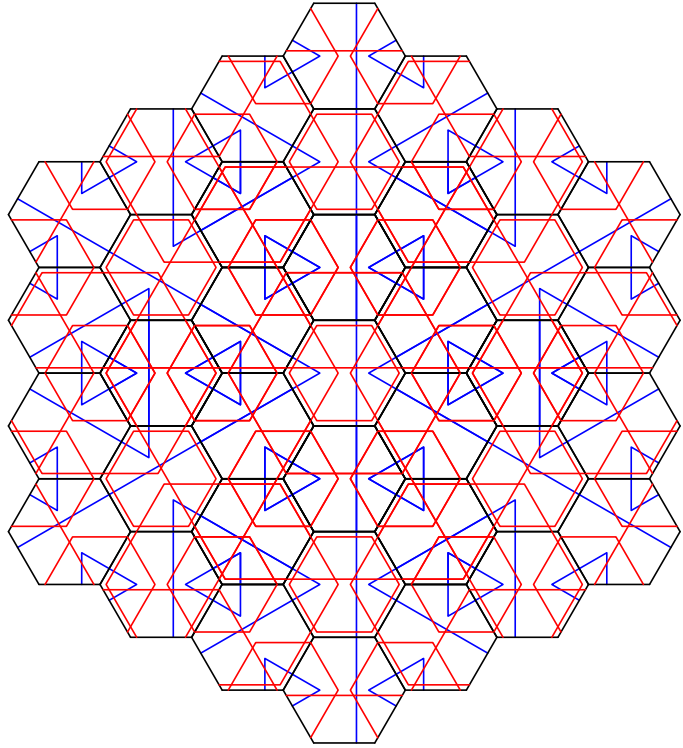}}
\caption{Patch of the $(1+\varepsilon + \varepsilon^{2})$-tiling in
  the variant with line decorations. In addition to the hierarchy of
  blue triangles from the half-hex tiling, there is now also a
  hierarchy of red hexagons. Hexagons of each size form a lattice
  periodic array, with lattices of different densities.}
\label{pen-linedeco}
\end{figure}

The hexagon of the arrowed half-hex tiling has still one remaining
line of mirror symmetry.  Both the
$(1+\varepsilon+\varepsilon^2)$-tiling and the Taylor tiling break
this remaining symmetry, alongside with a splitting into 7 subtypes of
hexagons, which differ in their behaviour under the inflation. For the
Taylor tiling, a pseudo inflation for fully asymmetric hexagons
follows easily from Figure~\ref{taysub}, which can then be converted
into a standard lattice inflation as before.

\begin{figure}[t]
\centerline{
\includegraphics[width=0.28\textwidth]{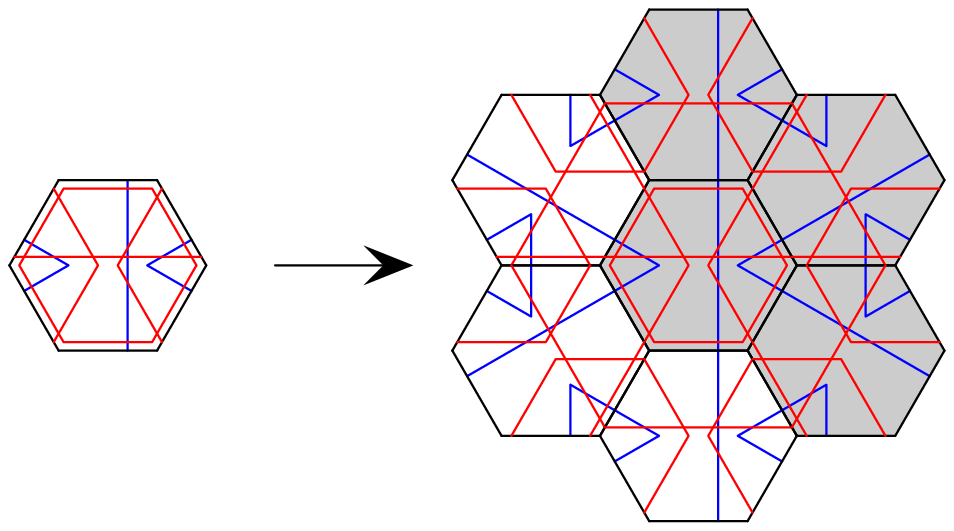}\hfill
\includegraphics[width=0.28\textwidth]{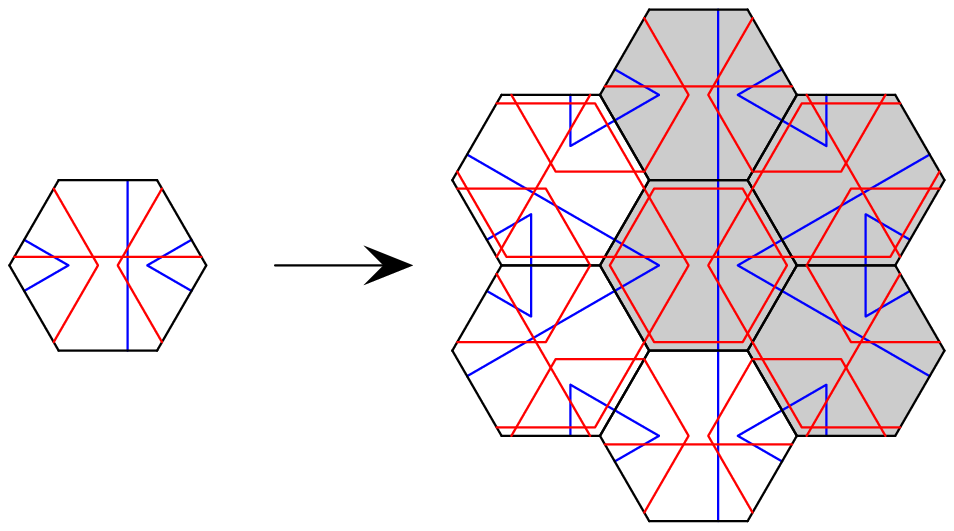}\hfill
\includegraphics[width=0.28\textwidth]{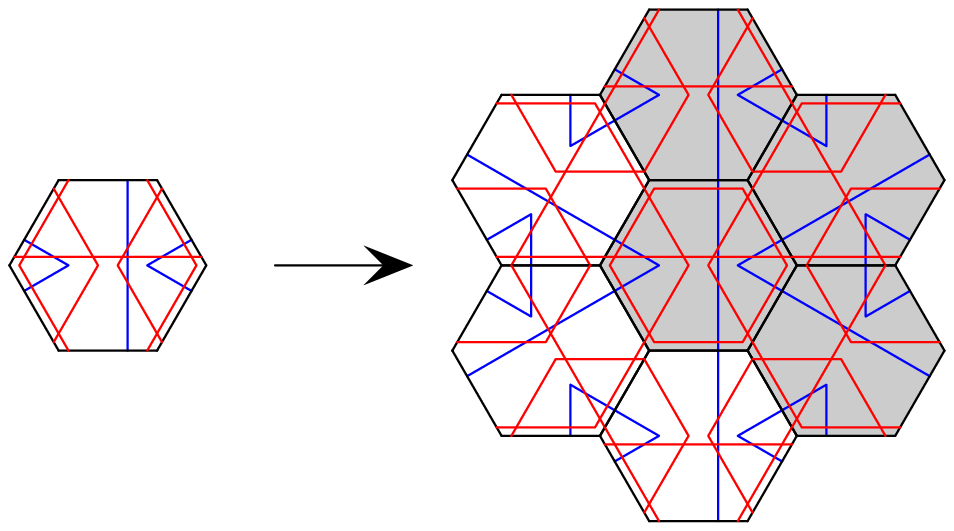}
}
\bigskip
\centerline{
\includegraphics[width=0.28\textwidth]{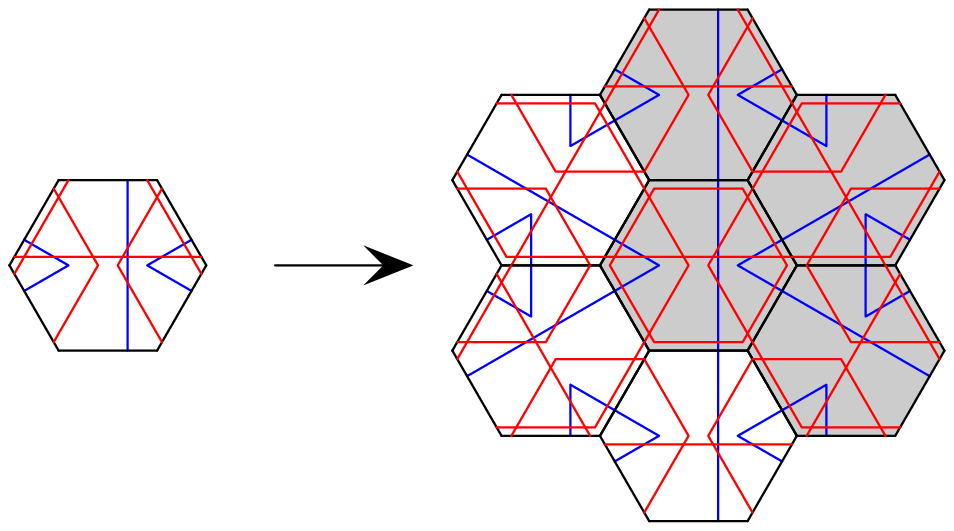}\hfil
\includegraphics[width=0.28\textwidth]{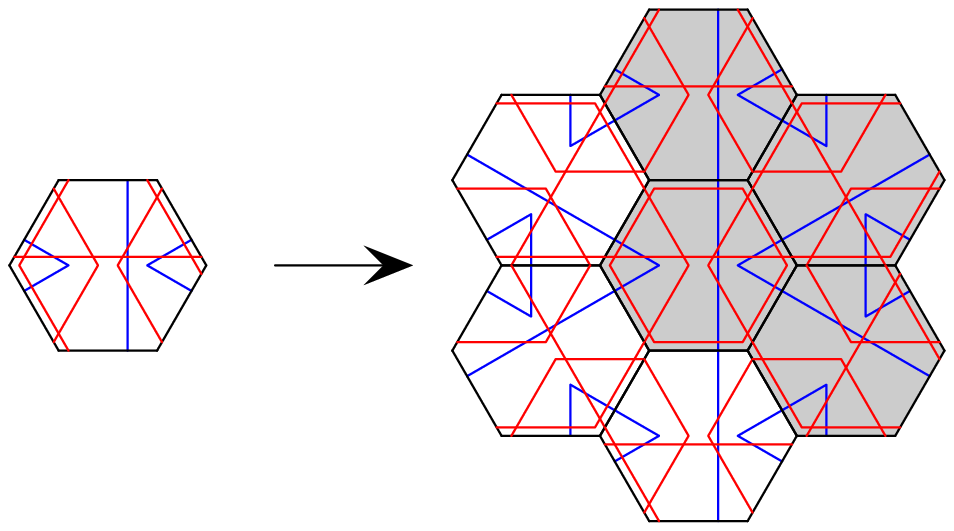}
}
\bigskip
\centerline{
\includegraphics[width=0.28\textwidth]{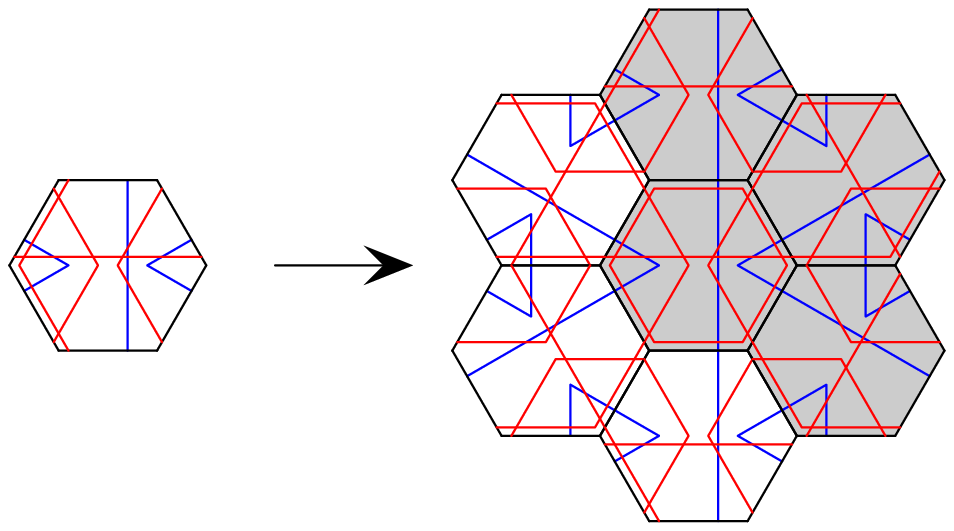}\hfil
\includegraphics[width=0.28\textwidth]{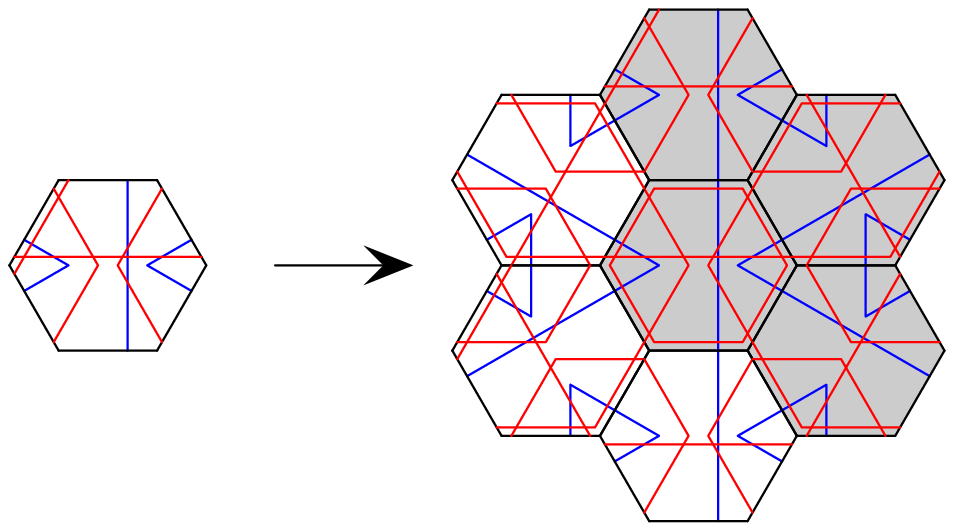}
}
\caption{Inflation rules for $(1+\varepsilon + \varepsilon^{2})$-tiling
in the variant with line decorations. If a hexagon is replaced by
7 hexagons, a pseudo inflation is obtained. Taking only the 4 shaded
hexagons results in a standard inflation rule (for one sector).}
\label{pen-linesubst}
\end{figure}

For the $(1+\varepsilon+\varepsilon^2)$-tiling, we replace the
$\varepsilon$- and $\varepsilon^2$-tiles by the equivalent line
decoration for the hexagons shown in Figure~\ref{pen-linedeco}, which
was introduced already in \cite{Pen97}. In addition to the blue lines
of the arrowed half-hex, which encode the arrows of the
$\varepsilon$-tiles along the hexagon edges, there is now also a
second set of red lines, which encode the information contained in the
$\varepsilon^2$-tiles and in the asymmetry across the
$\varepsilon$-tiles. The latter is represented by the piece of red
line parallel to the hexagon edge. The matching condition of the
$\varepsilon^2$-tiles is replaced by the requirement that red lines
must continue across tile edges. In the interior, each hexagon carries
an $X$-shaped pair of red line angles, and a red line belt, which
breaks the remaining mirror symmetry of the blue line decoration. The
inflation rule for the 7 types of decorated hexagons is shown in
Figure~\ref{pen-linesubst}.

{}From Figure~\ref{pen-linedeco}, it is clear that the red line
decoration breaks the horizontal mirror symmetry of the blue arrowed
half-hex decoration. Therefore, there will be three additional 1d
sub-solenoids $\mathbb{S}^{1}_{2}$ in three further directions, to
which the projection is 2-to-1 rather than 1-to-1. These 1d solenoids,
as well as those of the arrowed half-hex, all intersect in the single
point to which the now 12 fixed point tilings project.  Moreover, the
three `red' sub-solenoids intersect also in two further points at
corners of the hexagons, which are fixed points under the square of
the inflation. On these points, the projection to the solenoid is
6-to-1.

\begin{figure}[t]
\centerline{\includegraphics[width=0.4\textwidth]{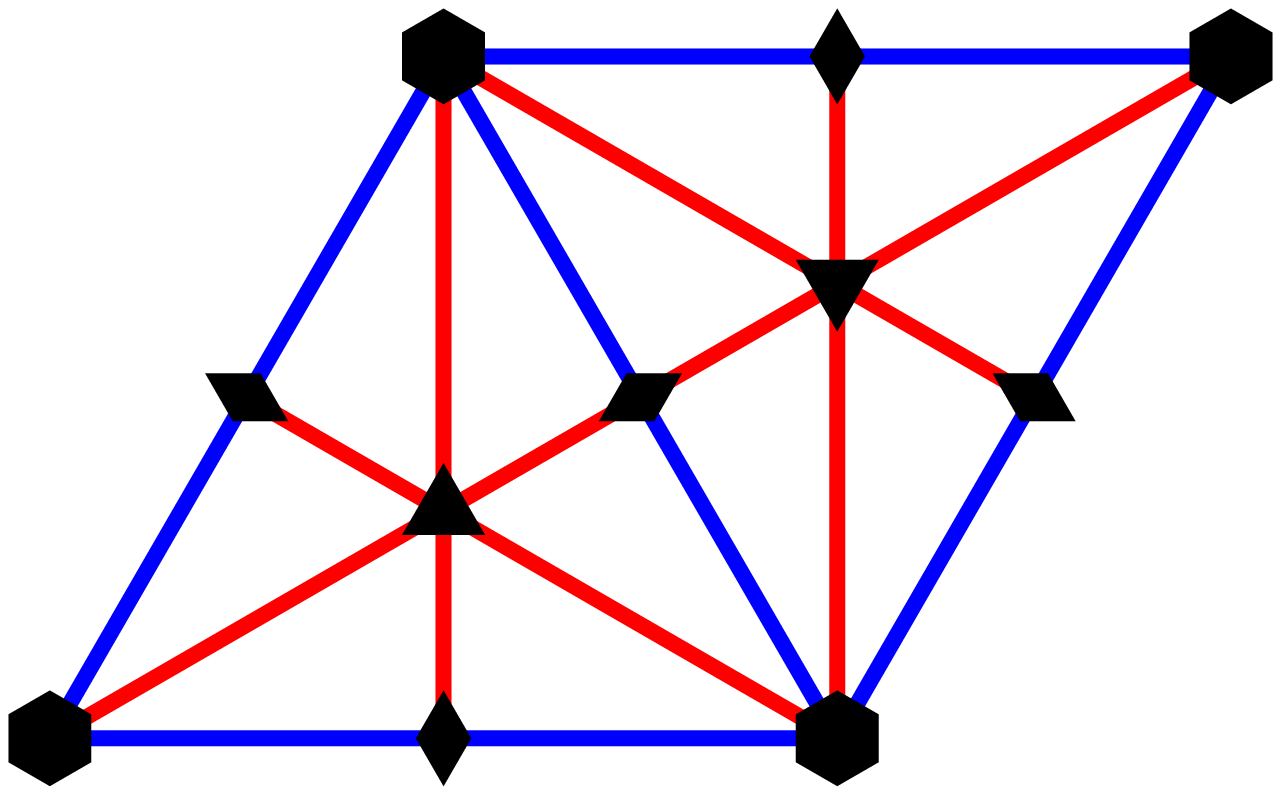}}
\caption{Schematic representation of the hulls of the different
  hexagon tilings via a (closed) toral slice of the underlying
  solenoid.  The projection fails to be 1-to-1 on the lines as
  follows.  For the \emph{half-hex}, the projection is 3-to-1
  precisely at the points marked by a hexagon, and 1-to-1 elsewhere.
  For the \emph{arrowed half-hex}, the 1d sub-solenoids where the
  projection is 2-to-1 are shown as blue lines.  At the hexagon
  points, three blue sub-solenoids intersect, and the projection is
  6-to-1.
  For the \emph{Taylor} and the \emph{Penrose} hexagon tilings, there
  are three additional sub-solenoids (red lines) along which the
  projection is 2-to-1. These intersect at two inequivalent points
  (marked by triangles) where the projection is 6-to-1 (at such
  points, three infinite order supertiles meet). The projection is
  12-to-1 at the hexagon points (which are centres of infinite-order
  supertiles).
  Note that, at the points marked by diamonds, not only two
  sub-solenoids intersect, but in fact all six. This is due to the
  dyadic structure of the hull. As these points have half-integer
  coordinates with respect to the hexagonal lattice, they are
  equivalent (in the solenoid) to the points marked by hexagons.}
\label{map}
\end{figure}

The set of points on the solenoid $\mathbb{S}^2_2$ to which the
projection of the Taylor tiling is not 1-to-1 is completely 
analogous. The corresponding ambiguities of the decoration,
starting with a fixed point tiling with full $D_6$ or $D_3$
symmetry, has been derived in \cite{ST11}. For an independent
approach, see \cite{LM}. A summary of the projection situation
is sketched and explained in Figure~\ref{map}.

The first benefit of having represented all four tilings
as lattice inflations is the following.

\begin{theorem}\label{pure-point}
  The half-hex, arrowed half-hex, $(1+\varepsilon + \varepsilon^{2})$
  and Taylor tilings are all model sets, and as such have pure point
  dynamical and diffraction spectrum.
\end{theorem}
\noindent\textsc{Sketch of Proof.}
As the half-hex and arrowed half-hex tilings are both factors of the 
other two tilings, we only need to consider the $(1+\varepsilon +
\varepsilon^{2})$- and Taylor tilings. Both have one `preferred' tile
type which occurs twice as often as the other tile types. For the
Taylor tiling, it is type $C$ (light blue in Figure~\ref{taysub}), and for
the $(1+\varepsilon + \varepsilon^{2})$-tiling it is the first type in
Figure~\ref{pen-linesubst}. In both cases, the centres of the preferred
tiles form a triangular sublattice of index $4$. Also the supertiles
(of any fixed order) of the preferred tile types form a periodic
array (disregarding orientation). Moreover, the preferred tiles are the 
seeds of the 12 fixed
point tilings, which differ from each other only along the 6 mirror
lines of the $D_6$ symmetry group. A sufficiently high order supertile
therefore contains tiles which are the same in all preferred
supertiles, independently of their orientation. Therefore, any Taylor
or $(1+\varepsilon + \varepsilon^{2})$-tiling contains a lattice-periodic 
subset of tiles. The results of \cite[Thm.~3]{LMold} then
imply that they must be model sets and hence have pure point 
dynamical and diffraction spectrum by \cite{LMS}.
\hfill$\square$\bigskip

Another benefit of dealing with a lattice inflation is that the
Anderson-Putnam approach \cite{AP98} to computing the \v{C}ech
cohomology of the hull becomes relatively easy to implement, since
tiles are represented as labelled lattice points whose environments
are easy to determine. Moreover, as the lattice inflation rules are
derived from overlapping pseudo inflation rules, they have the
property of forcing the border, which allows one to avoid the
complication of using collared tiles; see \cite{AP98} for details.

In this approach, from the local environments of the tiling a finite
approximant cell complex is constructed, whose points represent
cylinder sets of tilings. The full tiling space is then obtained as
the inverse limit space of the inflation acting on the approximant
cell complex. Correspondingly, the \v{C}ech cohomology of the hull is
the direct limit of the inflation action on the cohomology of the
approximant complex. For further details, we refer to
\cite{AP98,Sadun}; compare also \cite{BGG-gtm} for some explicitly
worked-out examples. For our four tiling spaces, the following results
are obtained.

\begin{theorem}
The four hexagon tiling spaces have the following \v{C}ech cohomologies:
\begin{fleqn}
\begin{alignat}{4}
&\text{Half-hex:} &\qquad
  H^2&=\mathbb{Z}[{\textstyle\frac14}]\oplus\mathbb{Z}^2, &\quad 
  H^1&=\mathbb{Z}[{\textstyle\frac12}]^2, &\quad H^0&=\mathbb{Z}; \notag \\
&\text{Arrowed half-hex:} &\qquad
  H^2&=\mathbb{Z}[{\textstyle\frac14}]
       \oplus\mathbb{Z}[{\textstyle\frac12}]^3, &\quad 
  H^1&=\mathbb{Z}[{\textstyle\frac12}]^2\oplus\mathbb{Z}, &\quad 
  H^0&=\mathbb{Z}; \notag \\
&\text{$1+\varepsilon+\varepsilon^2$:}  &\qquad
  H^2&=\mathbb{Z}[{\textstyle\frac14}]
       \oplus\mathbb{Z}[{\textstyle\frac12}]^6
       \oplus\mathbb{Z}^5\oplus\mathbb{Z}_3, &\quad 
  H^1&=\mathbb{Z}[{\textstyle\frac12}]^2\oplus\mathbb{Z}^2, &\quad  
  H^0&=\mathbb{Z}; \notag \\
&\text{Taylor:} &\qquad
  H^2&=\mathbb{Z}[{\textstyle\frac14}]
       \oplus\mathbb{Z}[{\textstyle\frac12}]^6\oplus\mathbb{Z}^7, &\quad 
  H^1&=\mathbb{Z}[{\textstyle\frac12}]^2\oplus\mathbb{Z}^4, &\quad 
  H^0&=\mathbb{Z}. \notag 
\end{alignat}
\end{fleqn}
\label{cohom}
In particular, the $H^{2}$-groups are distinct.
\hfill $\square$
\end{theorem}

Since the  $(1+\varepsilon + \varepsilon^{2})$- and the Taylor
tiling have different cohomology, it immediately follows that 
they cannot be in the same MLD class.

\begin{coro}
  The\/ $(1+\varepsilon + \varepsilon^{2})$ LI class and the
  Taylor LI class define distinct MLD classes.  \hfill $\square$
\end{coro}

As a side result of the cohomology calculation, the Artin-Mazur
dynamical zeta function of the inflation action on the hull
can also be obtained. It is defined as
\begin{equation}
 \zeta(z) = \exp{\left(\sum_{m=1}^\infty \frac{a_m}{m} z^m\right)}\notag
\end{equation}
where $a_m$ is the number of points in the hull invariant under an
$m$-fold inflation. We note that, if the hull consists of two (or
more) components for which the periodic points can be counted
separately, $a_m=a'_m+a''_m$, the corresponding partial zeta functions
have to be multiplied: $\zeta(z) = \zeta'(z)\cdot\zeta''(z)$. This
will turn out to be useful shortly. For the action of the mapping
$x \mapsto 2x$ on the dyadic solenoids $\mathbb{S}_{2}^{1}$ and
$\mathbb{S}_{2}^{2}$, the zeta function coincides with that of the
corresponding toral endomorphism, which are
\[
   \frac{1-z}{1-2z}
   \quad \text{and} \quad
   \frac{(1-2z)^{2}}{(1-z)(1-4z)}
\]
by an application of the results from \cite{BLP}. The corresponding
fixed point counts read $a_{m}^{(1)} = 2^{m} - 1$ and
$a_{m}^{(2)} = (2^{m} - 1)^{2}$, for $m\in\mathbb{N}$.

Anderson and Putnam \cite{AP98} have shown that the dynamical zeta
function can be computed from the action of the inflation on the
cochain groups of the approximant complex. Here, we rather express it
in terms of the action of the inflation on the rational cohomology
groups of the hull. If $A^{(m)}$ is the matrix of the inflation action
on the $m$-th rational cohomology group, the dynamical zeta function
is given by
\[
 \zeta(z) = \frac{\prod_{k\ \text{odd}}\det(1-zA^{(d-k)})}
                   {\prod_{k\ \text{even}}\det(1-zA^{(d-k)})}
          = \frac{\prod_{k\ \text{odd}}\prod_i(1-z\lambda_i^{(d-k)})}
                   {\prod_{k\ \text{even}}\prod_i(1-z\lambda_i^{(d-k)})}
\]
where the latter equality holds when the matrices $A^{(m)}$ are
diagonalisable with eigenvalues $\lambda_i^{(m)}$. The additional
terms in the expressions from \cite{AP98} cancel between numerator
and denominator.

{}From the cohomology and the eigenvalues $\lambda_i^{(m)}$, the dynamical
zeta function is easily obtained. For the half-hex, we get
\begin{equation}\label{zeta-hh}
  \zeta(z) = \frac{(1-2z)^2}{(1-z)(1-4z)(1-z)^2} =
  \frac{(1-2z)^2}{(1-z)(1-4z)} \cdot \frac{1}{(1-z)^2},
\end{equation}
which we have written as the product of the zeta function of the 2d
solenoid $\mathbb{S}_{2}^2$ and the generating function of two
additional fixed points. This is in line with our observation that the
projection to $\mathbb{S}_2^2$ is 1-to-1 except at one point, where it
is 3-to-1, wherefore there are two extra fixed points beyond the one
already present in $\mathbb{S}_2^2$.

For the arrowed half-hex, we find
\begin{equation}\label{zeta-ahh}
  \zeta(z) = \frac{(1-2z)^2(1-z)}{(1-z)(1-4z)(1-2z)^3} =
  \frac{(1-2z)^2}{(1-z)(1-4z)} \cdot {\left(\frac{1-z}{1-2z}\right)}^3 \cdot 
  \frac{1}{(1-z)^2},
\end{equation}
which is the product of the zeta functions of a 2d solenoid
$\mathbb{S}_2^2$, three 1d solenoids $\mathbb{S}^{1}_{2}$, and two
additional fixed points. As discussed above, we have three 1d
solenoids $\mathbb{S}^{1}_{2}$ where the projection to $\mathbb{S}_2^2$ is
2-to-1, hence the three extra copies of $\mathbb{S}^{1}_{2}$. The projection
from the fixed points of the arrowed half-hex tiling is 6-to-1 (they
form a $D_3$-orbit), so that there must be two extra fixed points
in the zeta function, in addition to those contained in the 4 solenoids.

Finally, both for the $(1+\varepsilon+\varepsilon^2)$- and the Taylor
tilings, we obtain
\begin{equation}\label{zeta-pt}
\begin{split} 
  \zeta(z) &= \frac{ (1-2z)^2 (1-z)^{2(+2)} }
      { (1-z) (1-4z) (1-2z)^6 (1-z)^{3(+2)} (1+z)^2 }  \\
 &= 
\frac{(1-2z)^2}{(1-z)(1-4z)} \cdot {\left(\frac{1-z}{1-2z}\right)}^6 \cdot 
\frac{1}{(1-z)^5} \cdot \frac{1}{(1-z^2)^2}, 
\end{split}
\end{equation}
which is the product of the zeta functions of a 2d solenoid
$\mathbb{S}_2^2$, six 1d solenoids $\mathbb{S}^{1}_{2}$, five additional
fixed points, and two extra 2-cycles. Even though the two tilings have
different cohomology, the additional terms in the zeta function of the
Taylor tiling (indicated by the extra exponents in parentheses) cancel
each other. As discussed above, both hulls contain six 1d solenoids
onto which the projection is 2-to-1.  The 12 fixed points, forming a
$D_6$-orbit, all project to the same point, where the six 1d
sub-solenoids intersect. In addition to the 7 fixed points in the
altogether 7 solenoids, there must hence be 5 further fixed points,
which indeed show up in the zeta function.  Finally, the three `red' 1d
solenoids intersect also at two types of corners of the hexagon
tiles. These points form 2-cycles under the inflation, and the
projection to them is 6-to-1 (two $D_3$-orbits).  Hence, two extra
2-cycles are present in the zeta function. The zeta functions derived
above confirm that our analysis of the set where the projection to
$\mathbb{S}_{2}^2$ fails to be 1-to-1 must have been complete.

\begin{coro}
  The Artin-Mazur zeta functions for the inflation action on the
  half-hex and on the arrowed half-hex hull are given by
  Eqs.~\eqref{zeta-hh} and \eqref{zeta-ahh}, while the $(1+\varepsilon
  +\varepsilon^{2})$-hull and the Taylor tiling hull have the same
  zeta function, as given by Eq.~\eqref{zeta-pt}.

  The corresponding fixed point counts, for $m\ge 1$, are given by\/
  $a_{m}^{(\mathrm{hh})} = (2^m - 1)^{2} + 2$, by\/
  $a_{m}^{(\mathrm{ahh})} = (2^m\! - 1)^{2} + 3 (2^{m}\! - 1) + 2$,
  and by\/ $a_{m}^{(\mathrm{P/T})} = (2^m\! - 1)^{2} + 6 (2^{m}\! - 1) + 5 +
  2 (1+(-1)^{m})$.  \hfill $\square$
\end{coro}

It is a rather amazing fact that the Penrose and the Taylor tiling,
despite defining distinct MLD classes, share the same dynamical
zeta function for the respective inflation action, and have a
projection to the 2d solenoid $\mathbb{S}_{2}^2$ with exactly the
same multiplicities.

\section{Outlook and open problems}

The two tiling spaces due to Penrose and to Taylor, which both are
model sets (see also \cite{LM}), show amazing similarities, though
they are certainly not MLD.  Whether there is a local derivation in
one direction, is still not fully clear, but unlikely. It is an open
problem where and what exactly is the difference between the two
tilings.

In the two parity patterns, beyond the percolation structure, one can
see islands of growing size in both tilings. Based on the inflation
structure, it is thus natural to conjecture that both classes of
parity patterns contain islands of unbounded size, though a proof does
not seem obvious. Also, the emergence via inflation series seems
slightly different in the two tilings.

The quest for a true monotile in the plane is not settled yet, because
the rules cannot be realised by nearest neighbour conditions, and
hence not by simple markings alone (unless one admits a prototile
version with disconnected parts).  There are other attempts to find an
example, for instance via polyominoes and related objects; compare
\cite{Rho05} and references therein.

An entirely different situation is met in $3$-space, where the famous
SCD prototile \cite{Danzer,BF} establishes a mechanism that is truly
three-dimensional. Indeed, the non-periodicity here is a result of a
screw axis with an incommensurate rotation, wherefore the repetitive
cases are aperiodic, but not strongly aperiodic; see \cite{BGnew} for
a discussion. Also, unlike the situation above, the local rules for
the SCD tile explicitly have to exclude the use of a reflected version,
which is perhaps not fully satisfactory either.

In summary, some progress was made in the quest for an aperiodic
monotile in recent years, but the search is certainly not over yet!

\section*{Acknowledgements}

We are grateful to Roger Penrose and Joan Taylor for important
comments and suggestions, to Robert Moody and Egon Schulte for helpful
discussions.  This work was supported by the German Research Council
(DFG), within the CRC 701. UG is grateful to the Fields Institute for
financial support.

\begin{small}

\end{small}

\end{document}